 \documentclass[11pt]{amsart}
 \usepackage{latexsym}
 \usepackage{graphics, epsfig}
 \newtheorem{theorem}{Theorem}[section]
 \newtheorem{proposition}{Proposition}[section]
 \newtheorem{lemma}{Lemma}[section]
 \newtheorem{corollary}{Corollary}[section]

 \newcommand{\rk}{\emph{Remarks}}

\begin{document}

 \title{THE GEOGRAPHY OF SPIN SYMPLECTIC $4$-MANIFOLDS}

 \author{Jongil Park}

 \address{Department of Mathematics, Konkuk University \\
          1 Hwayang-dong, Kwangjin-gu \\
          Seoul 143-701, Korea}

 \email{jipark@kkucc.konkuk.ac.kr}

 \thanks{This paper was supported by grant No. 1999-2-102-002-3
         from the interdisciplinary research program of  the KOSEF}

 \date{June 1, 2000 / Revised: July 23, 2001}



 \maketitle

\begin{abstract}
 In this paper we construct a family of simply connected spin
 non-complex symplectic $4$-manifolds which cover all but
 finitely many allowed lattice points ($\chi, {\mathbf c}$) lying
 in the region $0 \leq {\mathbf c} \leq  8.76\chi$.  Furthermore,
 as a corollary, we prove that there exist infinitely many
 exotic smooth \mbox{structures} on $(2n+1)(S^{2}\times S^{2})$
 for all $n$ large enough.
\end{abstract}

\section{Introduction}
\label{intro}

\markboth{JONGIL PARK}{THE GEOGRAPHY OF SPIN SYMPLECTIC
                       $4$-MANIFOLDS}

 One of the main topological problems in smooth $4$-manifolds is
 to classify such $4$-manifolds by using numerical invariants.
 Even though most classical invariants are not enough to
 distinguish smooth $4$-manifolds which are mutually homeomorphic,
 it is still of interest to know what combinations of topological
 invariants are realizable. Similarly, one may consider
 the same problems in both the category of symplectic $4$-manifolds and
 the category of complex surfaces. We call these geography problems.
 Let us describe some numerical invariants used in the geography question.
 Note that most  classical invariants for $4$-manifolds  are encoded
 by the intersection form. This  form is an integral unimodular symmetric
 bilinear pairing
 \[ Q_{X} : H_{2}(X; {\mathbf Z}) \times H_{2}(X; {\mathbf Z})
                                  \longrightarrow {\mathbf Z} \]
 obtained by representing homology classes as smoothly embedded oriented surfaces
 and counting intersections with signs. By Poincar\'{e} duality it is
 equivalently defined  as the cup product on the second cohomology classes.
 From $Q_{X}$ one can obtain the following  numerical invariants
 for a simply connected closed $4$-manifold $X$:

\begin{itemize}
 \item  $e(X) := b_{2}(X) + 2$,  called the Euler characteristic of $X$.
 \item  ${\bf \sigma }(X) := b_{2}^{+}(X) - b_{2}^{-}(X)$, called the signature
        of $X$.
 \item {\bf type} of $X$, which is {\bf even} if $Q_{X}(\alpha, \alpha) \in
        2{\mathbf Z}, \ \forall \alpha \in  H_{2}(X; {\mathbf Z})$, and {\bf
        odd} otherwise.
\end{itemize}

\hspace*{-1.8em} Now we define two numbers ${\mathbf \chi}(X)$ and ${\mathbf c}
  (X)$ for a simply connected closed  $4$-manifold  $X$ from these invariants:

\begin{itemize}
 \item $\chi(X) := \frac{\sigma(X) +e(X)}{4} = \frac{b_{2}^{+}(X) + 1}{2}$,
  which is an integer if there exists an almost complex structure on $X$, and
\end{itemize}

\begin{itemize}
 \item ${\mathbf c}(X) := 3\sigma(X) + 2e(X)$, which is equal to $c_{1}^{2}(X)\ $
      if there exists an almost complex structure on $X$.
\end{itemize}

\hspace*{-1.8em} Then a pair $(\chi(X), {\mathbf c}(X))$ of
 numbers determines the intersection form $Q_{X}$  up to type,
 so that a pair $(\chi(X), \mathbf{c}(X))$ determines the $4$-manifold
 $X$ up to homeomorphism and  type due to M. Freedman's classification
 theorem of simply connected topological $4$-manifolds.
 Thus it is natural to ask whether a pair $(\chi,\mathbf{c})$ of numbers
 determines a smooth (or symplectic) $4$-manifold in the category of
 smooth (or symplectic) $4$-manifolds.
 Explicitly, we may ask the following two questions:\\

\hspace*{-1.8em} {\bf Existence Problem} \hspace{2pt} {\em Which lattice points
  in the $(\chi, {\mathbf c})$-plane are realized as simply connected irreducible
  smooth $($or symplectic$)$ $4$-manifolds?}\\

\hspace*{-1.8em} {\bf Uniqueness Problem} \hspace{2pt} {\em If there exists a
 simply connected irreducible smooth $($or symplectic$)$ $4$-manifold
 corresponding to a given lattice point $(\chi, {\mathbf c})$, are there many
 distinct smooth $($or symplectic$)$ structures on it?}\\

 These are the geography problems of $4$-manifolds in which we are interested.
 Note that the questions above are meaningful in the category of
 irreducible $4$-manifolds. Otherwise, one can easily get an answer for those
 questions through the connected sum of other $4$-manifolds.
 Here we call a smooth $4$-manifold $X$ {\em irreducible}
 if it is not a connected sum of other smooth $4$-manifolds except for a
 homotopy $4$-sphere, i.e. if $X=X_{1}\sharp X_{2}$ implies that one of
 $X_{i}$ is a homotopy $4$-sphere.

  Despite the fact that it was a fundamental problem in the category of
 $4$-manifolds, the geography problem for smooth $4$-manifolds and for
 symplectic $4$-manifolds had long been unsolvable until gauge theory was
 introduced. Since the geography problem has almost the same answer
 even if we restrict ourselves to symplectic manifolds, from now on
 we only consider geography problems of symplectic $4$-manifolds.
 Since gauge theory was introduced by  S. Donaldson in 1982,
 this area has rapidly been developed. In particular, as an application
 of the Seiberg-Witten theory to symplectic $4$-manifolds,
 many  remarkable  results about symplectic $4$-manifolds have been obtained.
 For example,  C. Taubes  proved that every irreducible symplectic
 $4$-manifold $X$ with $b_{2}^{+}(X)  > 1$  has a non-zero Seiberg-Witten
 invariant and satisfies $\mathbf{c}(X)=c_{1}^{2}(X) \geq 0$ (\cite{t}).
 Besides, geography problems for irreducible symplectic $4$-manifolds
 have been extensively studied by topologists such as R. Fintushel,
 R. Gompf, \mbox{R. Stern,}  A. Stipsicz,  Z. Szab\'{o} and the author,
 so that most existence problems for irreducible symplectic $4$-manifolds
 are settled (\cite{fs1}, \cite{g}, \cite{p}, \cite{ps}, \cite{st}).
 But, when we restrict our concern to spin symplectic $4$-manifolds, it is
 still mysterious  which lattice points in the ($\chi, {\mathbf c}$)-plane
 are realized by simply connected spin symplectic $4$-manifolds.
  Note that every spin symplectic $4$-manifold is automatically
 {\em irreducible} because its Seiberg-Witten invariant is non-zero due to C.
 Taubes'  theorem above, so that, if it were not irreducible, it would be
 a blow up of another manifold  contradicting the spin hypothesis (\cite{d}).
 Since a spin structure on a manifold is a purely topological invariant
 (\cite{lm}), there is a topological restriction on
 \mbox{$\chi$ and ${\mathbf c}$}
 for spin $4$-manifolds. That is, every spin symplectic $4$-manifold $X$
 satisfies ${\mathbf c}(X) \equiv 8\chi(X)\!\!\! \pmod {16}$.
 Furthermore, if a $4$-manifold  $X$ has a \mbox{non-positive} signature,
 then it also satisfies $\mathbf{c}(X) \leq  8\chi(X)$.
 Hence the interesting questions are:
 Which lattice points in the ($\chi, {\mathbf c}$)-plane satisfying
 \[ 0 \leq  \mathbf{c} \leq 8\chi \ \ \mathrm{and} \ \
 \mathbf{c} \equiv 8\chi \!\! \pmod{16} \]
 are realized as $(\chi, c_{1}^{2})$ of simply connected spin symplectic
 $4$-manifolds, and then, if there exists such a spin symplectic $4$-manifold
 corresponding to a given lattice point $(\chi,\mathbf{c})$, are there many
 distinct symplectic structures on it?  One may also ask the same questions
 for lattice points lying in the positive signature region:
 ${\mathbf c} > 8\chi \ \ \mathrm{and} \ \  {\mathbf c} \equiv
    8\chi \!\! \pmod{16}.$

 In  this paper, we give an answer for these geography  problems.
 First  we  answer the existence question as follows:  Roughly,
 by taking a symplectic fiber sum along an embedded surfaces of
 self-intersection $0$ introduced by  R. Gompf (\cite{g}),
 we construct a family of simply connected spin symplectic $4$-manifolds
 which cover all `allowed' lattice points in the wedge between the
 elliptic line ($\mathbf{c}=0$) and a line  $\mathbf{c}= 2\chi -12$,
 which is parallel to the Noether line. And then, by taking a
 symplectic fiber sum repeatedly along a torus embedded in both such
 $4$-manifolds constructed above and an appropriate spin symplectic
 $4$-manifold with positive signature, we are able to construct desired
 simply connected spin symplectic $4$-manifolds.

 We also investigate the uniqueness question for symplectic
 $4$-manifolds. In general, since the Seiberg-Witten invariant of a smooth
 $4$-manifold with $b_{2}^{+} > 1$ is a diffeomorphic invariant, one can
 easily produce a family of smooth $4$-manifolds which are mutually
 homeomorphic, but which are not diffeomorphic by showing that they have
 different Seiberg-Witten invariants. One way to obtain such a family of smooth
 $4$-manifolds is to use {\em a logarithmic transformation} on a tubular
 neighborhood of a regular fiber lying in a cusp neighborhood (\cite{fs1},
 \cite{p}). Another way to obtain such a family  of smooth
 $4$-manifolds is to use {\em $0$-framed surgery}  on a knot $K$  embedded in
 $S^{3}$ (\cite{fs2}). Note that both techniques above can be
 performed  symplectically  in some cases.  For example, if $K$ is a fibered knot
 in $S^{3}$,  R. Fintushel and  R. Stern  constructed  a family of symplectic
 $4$-manifolds which are all homotopy $K3$ surfaces and which are pairwise
 non-diffeomorphic (\cite{fs2}). Since any two symplectic $4$-manifolds which are
 non-diffeomorphic  are automatically non-symplectomorphic,  their examples
 actually admit \mbox{infinitely} many distinct {\em symplectic} structures.
 We are going to use this technique to conclude that all simply connected spin
 symplectic $4$-manifolds obtained in the following main theorem admit
 infinitely many distinct symplectic structures.

\begin{theorem}
\label{main}
  There is a line ${\mathbf c} = f(\chi)$ with a slope $> 8.76$ in the
 $(\chi, {\mathbf c})$-plane  such that any allowed lattice point
 satisfying ${\mathbf c} \leq f(\chi)$ in the first quadrant
 can be realized as  $(\chi, c_{1}^{2})$ of  a simply connected spin
 non-complex symplectic $4$-manifold which admits infinitely many distinct
 symplectic structures. In particular, all allowed lattice points
 $(\chi, {\mathbf c})$ except finitely many lying in the region
 $0 \leq {\mathbf c} \leq 8.76\chi$ satisfy  ${\mathbf c} \leq f(\chi)$.
\end{theorem}

 \hspace*{-1.8em} {\em Remark.} \ B.D. Park and Z. Szab\'{o} got a
 similar result in~\cite{ps}. That is,
 they proved  that, for each even intersection form $Q$ satisfying
 $0 \leq {\mathbf c}(Q) :=  3\sigma(Q) + 2e(Q) < 8\chi(Q)$ and $b_{2}^{+}$
 odd, there is a simply connected irreducible spin  symplectic $4$-manifold
 $X$ with $Q_{X}=Q$. But they did not get spin symplectic $4$-manifolds
 with signature $0$ or positive signature. Theorem~\ref{main} above
 implies that, for every allowed lattice point $(\chi, \mathbf{c})$
 except finitely many lying in a non-negative signature region
 $8  \chi  \leq {\mathbf c} \leq  8.76\chi$, there exists a corresponding
 simply connected spin non-complex symplectic $4$-manifold
 which admits infinitely many distinct symplectic structures.\\

  Finally  we investigate exotic smooth structures on a connected sum
 $4$-manifold, denoted by $(2n+1)(S^{2}\times S^{2})$, of an odd number
 of copies of a smooth $4$-manifold $S^{2}\times S^{2}$.
 We say that a smooth $4$-manifold $X$ admits an {\em exotic} smooth structure
 if it has more than one distinct smooth structure, i.e. there exists a smooth
 $4$-manifold $X'$ which is homeomorphic to $X$, but not diffeomorphic to $X$.
 It has long been an interesting question whether a connected
 sum $4$-manifold $(2n+1)(S^{2}\times S^{2})$ admits an exotic smooth
 structure or not. In this paper, as a corollary of Theorem~\ref{main} above,
 we have the following remarkable result.

\begin{corollary}
\label{finale}
 There exists an integer $N$ such that, for all
 $\,n \geq N$, a connected sum $4$-manifold $(2n+1)(S^{2}\times S^{2})$
 admits infinitely many exotic smooth structures.
\end{corollary}

 This paper is organized as follows:  In Section~\ref{sec-2}
 we construct a family of simply connected spin non-complex
 symplectic $4$-manifolds shown in Theorem~\ref{main} above,
 and in Section~\ref{sec-3} we prove that all these manifolds
 satisfy desired properties stated in Theorem~\ref{main} above.\\

\section{Construction of spin symplectic $4$-manifolds}
\label{sec-2}
 We start this section by  briefly reviewing the geography problem
 of minimal complex surfaces. It is well known that every simply
 connected minimal complex surface is one of rational or ruled surfaces,
 elliptic surfaces and complex surfaces of general type.
 The geography problem for minimal complex surfaces of general type
 has been studied extensively by algebraic surface theorists.
 The following are some of main results:

\begin{itemize}
 \item  Chern-invariants $\chi =\frac{1}{12}(c_{1}^{2} + e)$ and
       $\mathbf{c} = c_{1}^{2}=3\sigma +2e$ of minimal complex surfaces
       of general type satisfy $\chi > 0$ and ${\mathbf c} > 0$.
       They also satisfy both Noether inequality and Bogomolov-Miyaoka-Yau
       inequality: $ 2\chi - 6 \leq \mathbf{c} \leq 9\chi $.
 \item  Every lattice point ($\chi, \mathbf{c}$) except finitely many lying
       in $2\chi -6 \leq \mathbf{c} \leq 8\chi$ is realized as
       ($\chi, c_{1}^{2}$) of a minimal complex surface which is
       a Lefschetz fibration. Furthermore, if it satisfies
       $2\chi -6 \leq {\mathbf c} \leq 8(\chi- \frac{9}{\sqrt[3]{12}}
       \chi ^{\frac{2}{3}})$, then it is realized as $(\chi, c_{1}^{2})$
       of a simply connected minimal complex surface which is a genus
       two Lefschetz fibration (\cite{per}).
 \item  By modifying G. Xiao's construction of non-spin
       complex surfaces which are hyperelliptic fibrations, U. Persson,
       C. Peters and G. Xiao constructed infinitely many simply connected
       spin complex surfaces with positive signature which are
       fibrations over ${\mathbf CP}^{1}$ even though these
       are not hyperelliptic fibrations. Roughly, these are constructed
       by using a triple sequence of double coverings
       $Y_{3} \stackrel{\pi_{3}} {\longrightarrow} Y_{2}
                \stackrel{\pi_{2}} {\longrightarrow} Y_{1}
                \stackrel{\pi_{1}} {\longrightarrow}
                       {\mathbf CP}^{1} \times {\mathbf CP}^{1}$
       branched at $B_{1}, \pi_{1}^{*}B_{2}$ and $(\pi_{1}\pi_{2})^{*}B_{3}$,
       respectively, where branch-loci
       $B_{1}, B_{2}$ and $B_{3}$ are chosen suitably in
       ${\mathbf CP}^{1} \times {\mathbf CP}^{1}$
       so that $Y_{3}$ is a simply connected spin singular complex
       \mbox{surface} with positive signature. After desingularization
       of $Y_{3}$ and blowing down disjoint exceptional curves,
       they obtained a Lefschetz fibration
       $Y \longrightarrow {\mathbf CP}^{1} \times {\mathbf CP}^{1}
       \stackrel{proj_{1}} {\longrightarrow} {\mathbf CP}^{1}$.
       Note that it is not hyperelliptic because its generic fiber is a
       ${\mathbf Z}_{2} \oplus {\mathbf Z}_{2} \oplus {\mathbf
       Z}_{2}$-cover over ${\mathbf CP}^{1}$
       (See \cite{pet}, \cite{ppx} for details).
\end{itemize}

  Now let us consider spin $4$-manifolds. Recall that an oriented manifold
 is called {\em  spin} if the  second  Stiefel-Whitney class $w_{2}$ of
 the manifold is zero. Equivalently, the intersection form of the manifold
 is of even type in the simply connected case. As we mentioned in the
 Introduction, the spin condition gives a topological restriction on
 $\chi(X):= \frac{\sigma(X)+e(X)}{4}$ and
 ${\mathbf c}(X):= 3\sigma(X) + 2e(X)$ for a spin $4$-manifold $X$.
 Explicitly,

\begin{lemma}
\label{lem-spin}
  If $X$ is a spin symplectic $4$-manifold with $b_{2}^{+} > 1$,
  then \mbox{it satisfies}
  \[ {\mathbf c}(X) = c_{1}^{2}(X) \geq  0  \ \ \ \mathrm{and} \ \ \
      {\mathbf c}(X) \equiv 8\chi (X)  \!\!\pmod{16}.\]
\end{lemma}

\begin{proof}
 \ The first inequality follows from  C. Taubes'  theorem (\cite{t}),
  and the second one follows from  the facts that
 ${\mathbf c}(X) = 3\sigma(X) + 2e(X) =\sigma(X) + 8\chi (X)$
 and $\sigma(X)$ is divisible by $16$ due to Rohlin's signature theorem.
\end{proof}

 Lemma~\ref{lem-spin} above tells us that only lattice
 points ($\chi, {\mathbf c}$) satisfying both ${\mathbf c} \geq 0$ and
 ${\mathbf c} \equiv 8\chi \!\!\! \pmod{16}$ can be possibly realized as
 $(\chi, c_{1}^{2})$ of spin symplectic $4$-manifolds - We call these
 `allowed' lattice points. But it is still not fully known which allowed
 lattice points in the ($\chi, \mathbf{c}$)-plane are covered by simply
 connected spin symplectic $4$-manifolds. Furthermore, there is another
 restriction on the geography of spin complex surfaces, that is, not all
 allowed lattice points lying in $2\chi -6\leq {\mathbf c} \leq 8\chi$
 are covered by spin complex surfaces. For example, U. Persson, C. Peters
 and G. Xiao proved:

\begin{theorem}[\cite{ppx}]
\label{thm-ppx}
  Let $X$ be a simply connected spin complex surface whose
  Chern-invariants satisfy $2\chi -6 \leq c_{1}^{2} < 3(\chi -5).$
  Then it has $c_{1}^{2} = 2(\chi -3)$ with $c_{1}^{2} = 8k$ and $k$ odd,
  or $c_{1}^{2} = \frac{8}{3}(\chi -4)$ with $\chi \equiv 0 \!\!\! \pmod{3}$.
\end{theorem}

  But the story is quite different in the category of symplectic $4$-manifolds.
 For example, we will show that all allowed lattice points lying in the region
 $2\chi -6 \leq c_{1}^{2} < 3(\chi -5)$ are also covered by simply connected
 spin symplectic $4$-manifolds. Before proving it, first we define a
 topological surgery, called {\em a fiber sum}, which is one of the main
 technical tools in our construction.\\

\hspace*{-1.8em} {\bf Definition} \
  Let $X$ and $Y$ be closed oriented smooth $4$-manifolds which contain
 a smoothly embedded surface $\Sigma$ with genus $g \geq 1$.
 Suppose $\Sigma$ represents a homology class of infinite order and of
 self-intersection $0$, so that there exists a tubular neighborhood,
 say $D^{2}\times \Sigma$, in both $X$ and $Y$.
 Let $X_{0}=\overline{X\setminus D^{2}\times \Sigma}$ and
 $Y_{0}=\overline{Y\setminus D^{2}\times \Sigma}$. Then by
 choosing an orientation-reversing, fiber-preserving diffeomorphism
 $\varphi: D^{2}\times \Sigma \longrightarrow D^{2}\times \Sigma$ and
 by gluing $X_{0}$ to $Y_{0}$ along their boundaries via the diffeomorphism
 \[\varphi |_{\partial D^{2}\times \Sigma}:
 \partial D^{2}\times \Sigma \longrightarrow  \partial D^{2}\times \Sigma \]
 we define a new oriented smooth $4$-manifold $X\sharp_{\Sigma}Y$,
 called the {\bf fiber sum} of $X$ and $Y$ along $\Sigma$.\\

  Note that there is an induced embedding of $\Sigma$ into
 $X\sharp_{\Sigma}Y$, well-defined up to isotopy, which represents
 a homology class of infinite order and of self-intersection $0$.
 Furthermore,  R. Gompf generalized {\em a fiber sum} technique for
 symplectic manifolds. That is, R. Gompf proved that if $\Sigma$ is a
 symplectic (or Lagrangian) surface of self-intersection $0$
 in symplectic $4$-manifolds $X$ and $Y$, then a fiber sum operation
 can be performed symplectically, so that the resulting manifold
 $X\sharp_{\Sigma}Y$ is symplectic. He also proved that a spin
 structure can be preserved under a fiber sum operation (\cite{g}).

 Here is a brief explanation for a spin structure on $X\sharp_{\Sigma}Y$:
 Suppose there are two spin structures on a tubular neighborhood
 $D^{2}\times \Sigma$ induced from $X$ and $Y$, respectively. Then
 the difference of two spin structures corresponds to an element
 of $H^{1}(D^{2}\times \Sigma ; {\mathbf Z_{2}})\cong
 H^{1}(\Sigma ; {\mathbf Z_{2}})$.
 Since the group of automorphisms, up to isotopy, of a bundle
 $D^{2}\times \Sigma \rightarrow \Sigma$ can be identified with
 $[\Sigma, S^1] \cong H^{1}(\Sigma ; {\mathbf Z})$ and since
 the coefficient homomorphism $H^{1}(\Sigma ; {\mathbf Z})\rightarrow
 H^{1}(\Sigma ; {\mathbf Z_{2}})$ is surjective, an element of
 $H^{1}(D^{2}\times \Sigma ; {\mathbf Z_{2}})$ corresponding to the
 difference of two spin structures is the image of a suitable bundle
 automorphism $\varphi$ under the coefficient homomorphism.
 Hence two spin structures on $D^{2}\times \Sigma$ can be identified
 by choosing a suitable bundle automorphism
 $\varphi : D^{2}\times \Sigma \rightarrow D^{2}\times \Sigma\,$ so that
 the spin structures on $X$ and $Y$ induce a spin structure on
 $X\sharp_{\Sigma}Y$.

\begin{lemma}[\cite{g}]
\label{lem-gompf}
  Suppose $X$ and $Y$ are spin symplectic $4$-manifolds containing
  a symplectic (or Lagrangian) surface $\Sigma$ of self-intersection $0$.
  Then the manifold $X\sharp_{\Sigma}Y$ is also spin and symplectic.
\end{lemma}

\begin{lemma}
\label{lem-sc}
 Suppose $X$ and $Y$ are simply connected smooth $4$-manifolds
 containing a tubular neighborhood $D^{2}\times \Sigma$ of an embedded
 surface $\Sigma$. If $S^{1}\cong \partial D^{2} \times \{pt\}$ bounds
 a disk in $X_{0}=\overline{X\setminus D^{2}\times \Sigma}$ or
 $Y_{0}=\overline{Y\setminus D^{2}\times \Sigma}$,
 then  the manifold $X\sharp_{\Sigma}Y$ is simply connected.
\end{lemma}

\begin{proof}
 \ This fact follows from Van Kampen's theorem:
 \[ \pi_{1}(X\sharp_{\Sigma}Y) = \frac{\pi_{1}(X_{0})\ast \pi_{1}(Y_{0})}{N} \]
 where `$\ast$' means a free product and $N$ is a  normal subgroup generated by
 $\pi_{1}(\partial D^{2} \times \Sigma)$.
\end{proof}

\begin{lemma}
\label{lem-cc}
 Let $X$ and $Y$ be closed simply connected smooth $4$-manifolds
 containing an embedded surface $\Sigma$ of self-intersection $0$
 and genus $g$. Then the manifold $X\sharp_{\Sigma}Y$ satisfies
 ${\mathbf c}(X\sharp_{\Sigma}Y)
 ={\mathbf c}(X) +{\mathbf c}(Y) +8(g-1)$ and
 $\chi(X\sharp_{\Sigma}Y) =\chi(X)+\chi(Y)+(g-1)$.
\end{lemma}

\begin{proof}
 \ This follows from the facts that
 $e(X\sharp_{\Sigma}Y) = e(X)+e(Y)+4(g-1)$ and
 $\sigma(X\sharp_{\Sigma}Y) = \sigma(X)+\sigma(Y)$.
\end{proof}

\hspace*{-1.8em} {\bf Brieskorn  manifold}. \ As the first step to
 construct desired spin symplectic $4$-manifolds, we consider a Brieskorn
 $4$-manifold defined as follows: For positive integers $p, q$ and  $r$,
 a Brieskorn $4$-manifold $B(p,q,r)$ is defined to be the Milnor fiber of
 the link of the isolated singularity of $z_{1}^{p}+z_{2}^{q}+z_{3}^{r}=0$
 in  ${\bf C}^{3}$. Then $B(p,q,r)$ is a smooth $4$-manifold with boundary
 and, if $p, q$ and $r$ are pairwise coprime, the boundary
 $\partial B(p,q,r) = \Sigma(p,q,r)$ is the corresponding Brieskorn homology
 $3$-sphere (\cite{b}). It is well known that each $B(p,q,r)$ has a natural
 compactification (by adding a complex curve at infinity) as a complete
 intersection in a weighted homogeneous space and that the singularities
 of this compactification can be resolved to obtain a simply connected
 algebraic surface $X(p,q,r)$ (\cite{eo}). The following fact will be
 used to construct our desired spin symplectic $4$-manifolds.

\begin{lemma}[\cite{fs1}]
\label{lem-br}
 If $(p,q,r)$ and $(p',q',r')$ are triples of positive integers with
 $p\leq p',$ $q \leq q'$  and $r\leq r'$, then $B(p,q,r) \subset B(p',q',r')$.
 Furthermore if $p\geq 2, \, q\geq 3$ and $r\geq 7$, then a Brieskorn
 $4$-manifold $B(p,q,r)$ contains a cusp neighborhood and also a $2$-sphere
 $S$ which intersects the cusp fiber transversely at a single point.
\end{lemma}

  A. Stipsicz pointed out that an embedded torus in the cusp
 neighborhood stated in Lemma~\ref{lem-br} above is a Lagrangian torus
 (\cite{st}). Here are several examples of Brieskorn $4$-manifolds
 which will serve as building blocks in our construction.\\

\hspace*{-1.8em} {\em Example 1.}
 \ It is known that a simply connected elliptic surface $E(n)$ with
 no multiple fibers and holomorphic Euler characteristic $n$ can be
 obtained as an algebraic surface $X(2,3,6n-1)$ which is diffeomorphic to
 \[  B(2,3,6n-1) \cup_{\Sigma(2,3,6n-1)} C(n) \]
 where $C(n)$, usually called {\em a Gompf nucleus}, is a neighborhood
 of a cusp fiber and a section which is an embedded $2$-sphere of
 self-intersection $-n$.\\

\begin{figure}[hbtp]
 \centerline{\epsfig{figure=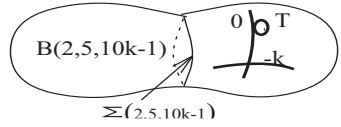, height=.7in, width=4.5in}}
 \caption{a Horikawa surface $H(4k\!-\!1)$}
 \label{Hori}
\end{figure}

\hspace*{-1.8em} {\em Example 2.}
 \ It is also known that a simply connected  Horikawa surface $H(4k-1)$,
 whose Chern-invariants are $(\chi, c_{1}^{2}) = (4k-1, 8k-8)$,
 can be obtained as an algebraic surface
  \[ X(2,5,10k - 1) = B(2,5,10k - 1)\cup_{\Sigma (2,5,10k - 1)} T(k)  \]
 where $T(k)$ is a neighborhood of a surface $T$ with genus $2$ and
 self-intersection $0$, obtained from $0$-framed surgery on the
 $(2,5)$-torus knot, and a section with self-intersection $-k$
 (See Fig.~\ref{Hori}).
 Note that a simply connected Horikawa surface $H(4k-1)$ is spin
 if and only if $k$ is even.\\

\begin{figure}[hbtp]
 \centerline{\epsfig{figure=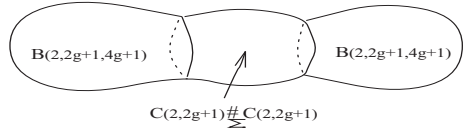, height=.7in, width=3.0in}}
 \caption{a symplectic $4$-manifold $Z$}
 \label{bries}
\end{figure}

\hspace*{-1.8em} {\em Example 3.}
 \ It is also known that a simply connected algebraic surface
 $X(2, 2g+1, 4g+1)$ is diffeomorphic to $B(2, 2g+1, 4g+1)
 \cup_{\Sigma (2, 2g+1, 4g+1)} T(2, 2g+1)$, where $T(2, 2g+1)$
 is a manifold obtained from $+1$-framed surgery on the
 $(2, 2g+1)$-torus knot so that it
 contains a surface $T$ of genus $g$ and self-intersection $0$.
 Let $X'=  X(2, 2g+1, 4g+1)\sharp \overline{\bf CP}^{2}$ be a
 manifold obtained by  blowing up at a point in $T$, so that
 $X'$ is diffeomorphic to
 \[ B(2, 2g+1, 4g+1)\cup_{\Sigma (2, 2g+1, 4g+1)} C(2, 2g+1) \]
 where $C(2, 2g+1)$ is a blow up of $T(2, 2g+1)$.
 Then $X'$ is a simply connected symplectic $4$-manifold which contains
 an embedded surface $\Sigma$ of genus $g$ and self-intersection $0$
 representing $T-e$.  Since $\Sigma$ is symplectically embedded,
 we get a simply connected symplectic $4$-manifold $Z:= X'\sharp_{\Sigma} X'$
 by taking a fiber sum of $X'$ with itself along $\Sigma$
 (See Fig.~\ref{bries}). Furthermore, the manifold $Z$ is spin because
 each part in
 Fig.~\ref{bries} is  spin and the boundary of each part in
 Fig.~\ref{bries} is a homology $3$-sphere $\Sigma (2, 2g+1, 4g+1)$.
 Note that the manifold $Z$ has Chern-invariants
 $\chi(Z)=2g^2-g+1$ and ${\mathbf c}(Z)=8g^2-16g+8$.\\

 Now we are ready to get our main theorem. First we prove the
 following proposition.

\begin{proposition}
\label{main-pro}
  Every allowed lattice point  $(\chi, {\mathbf c})$ in the  wedge between
 the elliptic line $({\mathbf c}=0)$ and a line ${\mathbf c} =2\chi-12$
 is realized as $(\chi, c_{1}^{2})$ of a simply
 connected spin non-complex symplectic $4$-manifold.
\end{proposition}

\begin{proof}
 \ First we claim that every  lattice point  $(\chi, {\mathbf c})$ satisfying
 both ${\mathbf c} \equiv 8\chi \pmod{16}$ and $\chi$ odd  on the Noether line
 ${\mathbf c} =2\chi - 6$ is realized as $(\chi, c_{1}^{2})$ of a  simply
 connected spin Horikawa surface  $H(8k'-1)$.  This follows  from the fact that
 $H(8k'-1)$  has Chern-invariants $(\chi, c_{1}^{2}) = (8k'-1, 16k'-8)$.
 And then, since a Brieskorn $4$-manifold  $B(p,q,r)$ with $(p,q,r) \geq (2,3,7)$
 contains an embedded Lagrangian torus $f$ (\cite{st}), we can prove by taking
 a symplectic fiber sum along the Lagrangian torus $f$ embedded in both
 $B(2,5,20k'-1) \subset H(8k'-1)$ and $B(2,3,12n-1)\subset E(2n)$
 that all allowed lattice points lying in the region
 $\{(\chi, {\mathbf c}) \  | \ 0 \leq {\mathbf c} \leq 2\chi - 6 \ \
  \mathrm{and} \ \ {\mathbf c} \equiv 8 \!\!\! \pmod{16} \}$
  \mbox{are covered by} a family of the following  spin symplectic
  $4$-manifolds
 \[ \{ H(8k'-1) \sharp_{f} E(2n) \ | \  k' \ \mathrm{and}
    \ n \ \mathrm{are \ positive \ integers} \},\]
 and all allowed lattice points in
  $\{(\chi, {\mathbf c}) \, | \, 0 \! \leq \! {\mathbf c}\! \leq \!2\chi
 \!-\! 12 \, \, \mathrm{and} \, \, {\mathbf c} \equiv 0 \!\!\! \pmod{16} \}$
 are also covered by the following  spin symplectic $4$-manifolds
 \[ \{ H(7) \sharp_{f}  H(8k' -1) \sharp_{f} E(2n) \ | \  k'
    \ \mathrm{and} \  n \ \mathrm{are \ positive \ integers}  \}.\]
 Note that all these manifolds are simply connected
 due to Lemma~\ref{lem-sc} and cannot admit a complex structure because they
 lie below the Noether line. Furthermore, we claim that every allowed lattice
 point $(\chi,{\mathbf c})=(2n,0)$ on the elliptic line $({\mathbf c}=0)$
 is also covered by a spin non-complex symplectic $4$-manifold $X(2n)$
 which is homeomorphic to, but not diffeomorphic to $E(2n)$ (See below).
 Thus we arrive at Proposition~\ref{main-pro} by combining all these results. \\

\hspace*{-1.8em} {\em Construction of $X(2n)$} :  Since there are
 cusp neighborhoods $N \subset B(2,3,11)$  and  $N' \subset B(2,3,12n-13)$,
 we construct a spin symplectic $4$-manifold $X(2n):= E(2)\sharp_{f}E(2n-2)$
 by taking a fiber sum along a Lagrangian torus $f$ embedded in both
 $N \subset B(2,3,11) \subset E(2)$ and $N' \subset B(2,3,12n-13)
 \subset E(2n-2)$. That is, $X(2n)$ is diffeomorphic to
 \[ C(2) \cup_{\Sigma{(2,3,11)}} B(2,3,11)\sharp_{f} B(2,3,12n-13)
    \cup_{\Sigma{(2,3,12n-13)}} C(2n-2) \]
 where $C(2)$ and $C(2n-2)$ are  nuclei, i.e. neighborhoods of
 a cusp fiber and a section in $E(2)$ and $E(2n-2)$, respectively.
 Then it is easy to check that both $X(2n)$ and $E(2n)$ have the same
 $(\chi,c_{1}^{2}) =(2n,0)$, but they are not diffeomorphic because
 they have different Seiberg-Witten invariants (Refer to Theorem~\ref{mms}
 of this paper and to \cite{fs1.5} for the Seiberg-Witten invariant
 of $X(2n)$ and $E(2n)$, respectively). Furthermore, a manifold  $X(2n)$
 cannot admit a complex structure because
 its Seiberg-Witten invariant is different from that of any
 elliptic surface (\cite{fs1.5}).
\end{proof}

\hspace*{-1.8em}  Let $\Omega$ be the set of all allowed lattice points
 lying in the wedge between the elliptic line $({\mathbf c}=0)$ and
 a line ${\mathbf c} =2\chi-12$. Then we have:

\begin{corollary}
\label{main-cor}
  If $X$ is a simply connected spin symplectic $4$-manifold which contains
 a Lagrangian torus $f$ in a cusp neighborhood, then each lattice point
 lying in $(\chi(X), \mathbf{c}(X))+\Omega$ is realized as $(\chi, c_{1}^{2})$
 of a simply connected spin symplectic $4$-manifold.
\end{corollary}

\begin{proof}
 \ By taking a symplectic fiber sum along a Lagrangian torus  embedded
 in a cusp neighborhood both in $X$ and in spin symplectic $4$-manifolds
 constructed in Proposition~\ref{main-pro} above, we get a family of
 desired simply connected spin symplectic $4$-manifolds.
\end{proof}

  Finally, by taking a fiber sum along a symplectic surface
 $\Sigma$ embedded in both a spin symplectic $4$-manifold $Z$ constructed
 in Example 3 above and an appropriate spin symplectic $4$-manifold
 of positive signature, we construct symplectic $4$-manifolds which appear
 in the main theorem of this paper:

\begin{theorem}
\label{main-1}
 There is a line ${\mathbf c} = f(\chi)$ with a slope $> 8.76$ in the
 $(\chi, {\mathbf c})$-plane such that any allowed lattice point
 $(\chi, {\mathbf c})$ satisfying ${\mathbf c} \leq f(\chi)$ in the first
 quadrant can be realized as $(\chi, c_{1}^{2})$ of  a simply connected
 spin  symplectic $4$-manifold. In particular, all allowed lattice points
 $(\chi, {\mathbf c})$ except finitely many lying in the region
 $0 \leq {\mathbf c} \leq 8.76\chi$ satisfy ${\mathbf c} \leq f(\chi)$.
\end{theorem}

\begin{proof}
 \ First let us choose a simply connected spin complex surface, say $Y$,
 of positive signature which contains a holomorphic curve $\Sigma_{g}$
 of self-intersection $0$ and genus $g$, and which also contains an embedded
 $2$-sphere $S$ intersecting $\Sigma_{g}$ transversely at one point.
 For example, we can choose such a complex surface $Y$ that U. Persson,
 C. Peters  and G. Xiao constructed in \cite{ppx}. Here is a brief
 construction of $Y$: Keeping in mind U. Persson, C. Peters and
 G. Xiao's construction noted at the beginning of Section~\ref{sec-2},
 let us consider an icosahedral group
 $I_{60}$ of order $60$ acting on ${\mathbf CP}^{1}$ and, for each
 integer $x> 0$, let $\Gamma(x)$ be a pull back of $\Gamma$ under
 a map ${\mathbf CP}^{1} \times {\mathbf CP}^{1} \longrightarrow
 {\mathbf CP}^{1} \times {\mathbf CP}^{1}$ induced from a standard
 $x$-covering of ${\mathbf CP}^{1}$ to itself,
 where $\Gamma = \bigcup_{g\in I_{60}} \{(p,g\cdot p)\in
 {\mathbf CP}^{1} \times {\mathbf CP}^{1} |\,
 p\in {\mathbf CP}^{1}\}$ is the union of all graphs obtained by
 the action of $I_{60}$ on ${\mathbf CP}^{1}$.
 Then $\Gamma(x)$ is a singular curve of bidegree $(60x,60x)$ in
 ${\mathbf CP}^{1} \times {\mathbf CP}^{1}$. Now choose
 branch-loci $B_{1},B_{2}$ and $B_{3}= \Gamma(x)\cup (15x,15x)$
 as singular curves of bidegree $(62x,*),(*,62x)$ and $(75x,75x)$ in
 ${\mathbf CP}^{1} \times {\mathbf CP}^{1}$, respectively.
 Then a triple sequence of double coverings
 $Y_{3} \stackrel{\pi_{3}} {\longrightarrow} Y_{2} \stackrel{\pi_{2}}
  {\longrightarrow} Y_{1} \stackrel{\pi_{1}} {\longrightarrow}
  {\mathbf CP}^{1} \times {\mathbf CP}^{1}$  branched at
  $B_{1}, \pi_{1}^{*}B_{2}$ and $(\pi_{1}\pi_{2})^{*}B_{3}$
  with minor modifications, respectively, induces
  a Lefschetz fibration $Y \longrightarrow {\mathbf CP}^{1}$
  whose Chern-invariants are asymptotically
  $\chi(Y) \approx 6857x^2$ and ${\mathbf c}(Y)= c_{1}^2(Y)
  \approx 60068x^2$ for sufficiently large $x$.

 Next take any simply connected spin symplectic $4$-manifold $Z$ which
 contains the same symplectically embedded surface $\Sigma_{g}$ of
 self-intersection $0$ as $Y$ has and which also contains a symplectic
 (or Lagrangian) torus $f$ in a cusp neighborhood $N$  satisfying
 $N\cap \Sigma_{g}=\phi$. For example, we choose such a $4$-manifold
 $Z$ as constructed in Example 3 above. And then, we take a symplectic
 fiber sum of $Y$ and $Z$ along the embedded surface $\Sigma_{g}$ to get
 a spin symplectic $4$-manifold $Y\sharp_{\Sigma_{g}} Z$.
 Since $[\mathbf{c}(Y)-8\chi(Y)] > 0$ and
 ${\mathbf c}(Y\sharp_{\Sigma_{g}} Z)-8\chi(Y\sharp_{\Sigma_{g}} Z) =
    [{\mathbf c}(Y)-8 \chi(Y)] +[{\mathbf c}(Z)-8\chi(Z)]$,
 there exists an integer $k>0$ such that $X := \overbrace{Y\sharp_{\Sigma_{g}}
 \cdots \sharp_{\Sigma_{g}} Y}^{k} \sharp_{\Sigma_{g}} Z$  has a positive
 signature. Note that the manifold $X$ has Chern-invariants
 $\chi(X)=k\chi(Y)+\chi(Z)+k(g-1)$ and
 ${\bf c}(X)=k{\bf c}(Y)+{\bf c}(Z)+8k(g-1)$,
 so its ratio is approximately
 \[ \frac{{\bf c}(X)}{\chi(X)} =
    \frac{k{\bf c}(Y)+{\bf c}(Z)+8k(g-1)}{k\chi(Y)+\chi(Z)+k(g-1)}
    \approx \frac{{\bf c}(Y)}{\chi(Y)} \approx \frac{60068x^2}{6857x^2}
    = 8.76009\cdots \]
 for sufficiently large integers $k$ and $x$.
 We fix such large integers $k$ and $x$ making
 $\frac{{\bf c}(X)}{\chi(X)} > 8.76$.
 In addition, since the complex surface $Y$ contains a
 $2$-sphere $S={\mathbf CP}^{1}$ intersecting
 a generic fiber $\Sigma_{g}$ transversely at one point,
 Lemma~\ref{lem-sc} implies that the spin symplectic $4$-manifold $X$
 constructed above is also simply connected.
 Let $\Omega$ be the set of all allowed lattice points lying
 in the wedge between the elliptic line (${\mathbf c}= 0$) and a line
 \mbox{${\mathbf c} = 2\chi-12$}.
 Then Corollary~\ref{main-cor} implies that every allowed lattice point
 lying in the region $(\chi(X), {\mathbf c}(X))+\Omega$
 is covered by a simply connected spin symplectic $4$-manifold.
 Furthermore, it is also true for $X\sharp_{f} X, \ X\sharp_{f}
 X\sharp_{f} X,\ldots \ $
 Hence, if we define a line ${\mathbf c}=f(\chi)$ by
 \[ f(\chi) = {\mathbf c}(X)/\chi(X) \cdot [\chi -{\mathbf c}(X)/2 - 6]
 + {\mathbf c}(X)\, , \]
 then each allowed lattice point $(\chi, {\mathbf c})$
 satisfying ${\mathbf c} \leq f(\chi)$ in the first quadrant is realized
 as $(\chi, c_{1}^{2})$ of a simply connected spin symplectic $4$-manifold
 $W:=\overbrace{X\sharp_{f} X\sharp_{f} \cdots \sharp_{f} X}^{m}
 \sharp_{f} V$ for some $m \in {\mathbf Z}$ and the manifold
 $V=H(8k' -1) \sharp_{f} E(2n)$
 (or $V=H(7)\sharp_{f} H(8k' -1) \sharp_{f} E(2n)$)
 constructed in the proof of Proposition~\ref{main-pro} above.
 Thus the theorem follows from the fact that the slope of
 $f(\chi) = {\mathbf c}(X)/\chi(X)$ is greater than $8.76$.
\end{proof}

\hspace*{-1.8em} \rk\
 1. Note that every spin symplectic $4$-manifold with signature $0$
 obtained in Theorem~\ref{main-1} above is homeomorphic to a connected
 sum $4$-manifold $(2n+1)(S^{2}\times S^{2})$ for some $n$ due to
 M. Freedman's classification theorem of simply connected closed topological
 $4$-manifolds. We will prove in Section~\ref{sec-3} that all those connected
 sum $4$-manifolds have infinitely many exotic smooth structures
 (See Corollary~\ref{main-2}).\\
 2. Theorem~\ref{main-1} above also provides a partial answer on the
 existence question of spin symplectic $4$-manifolds having a positive
 signature. That is, all but finitely many allowed lattice points
 $(\chi, {\mathbf c})$ lying in the region \mbox{$8\chi < {\mathbf c}
 \leq  8.76\chi$} are realized as $(\chi, c_{1}^{2})$ of simply connected
 spin symplectic $4$-manifolds of positive signature.\\
 3. There are still  infinitely  many  lattice points
 $(\chi,{\mathbf c})$
 lying in the region $f(\chi)< {\mathbf c} < 9\chi $ which are realized
 as $(\chi, c_{1}^{2})$ of simply connected spin non-complex
 symplectic $4$-manifolds. Furthermore, we do not claim that the line
 ${\mathbf c} = f(\chi)$ constructed in the proof above is the best
 choice.\\
 \\

\section{Proof of Theorem~\ref{main}}
\label{sec-3}
  In this section we complete the proof of Theorem~\ref{main} stated
 in the Introduction. There are two things left to be done - the uniqueness
 question of symplectic structures and the existence question of a complex
 structure on those symplectic $4$-manifolds that appear in
 Theorem~\ref{main}. First we prove:
\\

\hspace*{-1.8em} {\bf Claim 1}. {\em  Every simply connected spin
 symplectic $4$-manifold  in Theorem~\ref{main-1} admits infinitely many
 distinct symplectic structures.} \\

  As a background for the proof of Claim 1, let us introduce R. Fintushel and
 R. Stern's technique (See \cite{fs2} for details): Suppose $K$ is a fibered
 knot in $S^{3}$ with a punctured  surface $\Sigma_{g}^{\circ}$ of genus $g$
 as a fiber. Let $M_{K}$ be a $3$-manifold obtained by performing $0$-framed
 surgery on $K$, and let $m$ be a meridional circle to $K$.
 Then the $3$-manifold $M_{K}$ can be considered as a fiber bundle over
 a circle $m$ with a closed Riemann surface $\Sigma_{g}$ as a fiber, and there
 is a smoothly embedded torus $T_{m}:= m \times
 S^{1}$ of self-intersection $0$ in $M_{K} \times S^{1}$. Thus $M_{K} \times
 S^{1}$ is a fiber bundle over $S^{1}\times S^{1}$ with $\Sigma_{g}$ as a fiber
 and with $T_{m}= m \times S^{1}$ as a section. A theorem of Thurston states
 that such a $4$-manifold $M_{K} \times S^{1}$ has a symplectic structure with
 a symplectic section $T_{m}$ (\cite{th}). Thus, if $X$ is a symplectic
 $4$-manifold with a symplectically embedded torus $T$ of
 \mbox{self-intersection $0$,}  then the fiber sum $4$-manifold
 $X_{K} := X\sharp_{T=T_{m}} (M_{K} \times S^{1})$, obtained by taking
 a fiber sum along $T=T_{m}$, is symplectic. R. Fintushel and R. Stern
 proved that $X_{K}$ is homotopy equivalent to $X$ under a mild condition
 on $X$, and they also computed the Seiberg-Witten invariant of $X_{K}$.
 Explicitly,

\begin{theorem}[\cite{fs2}]
\label{fs-2}
  Suppose $X$ is a simply connected symplectic $4$-manifold which contains
 a symplectically embedded torus $T$ of self-intersection $0$ in a cusp
 neighborhood with $\pi_{1}(X\setminus T) = 1$ and representing a non-trivial
 homology class $[T]$. If $K$ is a fibered knot, then $X_{K}$ is
 a symplectic $4$-manifold which is homeomorphic to $X$ and whose
 Seiberg-Witten invariant is
 \[ SW_{X_{K}} = SW_{X} \cdot \Delta_{K}(t) \]
 where $\Delta_{K}(t)$ is the Alexander polynomial of $K$ and $t= \exp(2[T])$.
\end{theorem}

\hspace*{-1.8em} {\em Proof of Claim 1.} \ Since any two
 symplectic $4$-manifolds which have different Seiberg-Witten invariants
 are automatically non-symplectomorphic, we will apply Theorem~\ref{fs-2}
 to prove Claim 1. For this, it suffices to show that every
 simply connected spin symplectic $4$-manifold, say $W$, in
 Theorem~\ref{main-1} contains a symplectically embedded torus $T$ of
 self-intersection $0$ with $\pi_{1}(W\setminus T)= 1$ and representing
 a non-trivial homology class $[T]$.
 Note that each symplectic $4$-manifold shown in Theorem~\ref{main-1}
 is of the form
\begin{eqnarray*}
  W &=& \overbrace{X\sharp_{f} X\sharp_{f} \cdots \sharp_{f} X}^{m} \sharp_{f}
           H(8k' -1) \sharp_{f} E(2n)\, \, \, \, \, \, \,\, \, \mathrm{or}  \\
    &=& \overbrace{X\sharp_{f} X\sharp_{f} \cdots \sharp_{f} X}^{m}
           \sharp_{f} H(7)\sharp_{f} H(8k' -1) \sharp_{f} E(2n)
\end{eqnarray*}

\hspace*{-1.75em} for some positive integers $m, n$ and $k' \in {\mathbf Z} $,
 and an elliptic surface $E(2n)$ obviously contains such a symplectically
 embedded torus $T$ in a nucleus $C(2n)$ (Refer to Example 1 in
 Section~\ref{sec-2}). Thus Claim 1 follows from Theorem~\ref{fs-2} above.
 $\ \ \ \Box$  \\

  Now, as an application of Claim 1, we get an answer for the problem
 regarding how many exotic smooth structures exist on a connected sum of
 an odd number of copies of a $4$-manifold $S^{2}\times S^{2}$, denoted by
 $(2n+1)(S^{2}\times S^{2})$.

\begin{corollary}
\label{main-2}
 There exists an integer $N$ such that, for all
 $\,n \geq N$, a connected sum $4$-manifold $(2n+1)(S^{2}\times S^{2})$
 admits infinitely many exotic smooth structures.
\end{corollary}

\begin{proof}
 \ Note that each allowed lattice point $(\chi, {\mathbf c})$ on the signature
 $0$ line is realized by a simply connected spin smooth $4$-manifold
 $(2n+1)(S^{2}\times S^{2})$ for some $n$ and any simply connected spin
 smooth $4$-manifold having numerical invariants $(\chi, {\mathbf c}) =
 (n+1,8n+8)$ is homeomorphic to a spin $4$-manifold $(2n+1)(S^{2}\times S^{2})$
 due to M. Freedman's classification. Hence, by applying these facts to
 simply connected spin symplectic $4$-manifolds of signature $0$ shown
 in Theorem~\ref{main-1}, we get the desired result.
\end{proof}

\hspace*{-1.8em} \rk\
 1.  A connected sum $4$-manifold $(2n+1)(S^{2}\times S^{2})$
 has a natural smooth structure which is compatible with a standard smooth
 structure on each $4$-manifold $S^{2}\times S^{2}$, and its Seiberg-Witten
 invariant is zero with respect to this natural smooth structure.
 Hence we realize that any exotic smooth structure on
 $(2n+1)(S^{2}\times S^{2})$ obtained in Corollary~\ref{main-2}
 above is strikingly different from the standard one.\\
 2.  The constant $N$ in Corollary~\ref{main-2} can be any number
  $\chi$ satisfying  $f(\chi) \geq  8\chi$, where ${\mathbf c} = f(\chi)$
  is a line stated in Theorem~\ref{main-1}. Practically, $N$ will be
  a huge number because we have to choose sufficiently large
  integers $k$ and $x$ in order to get a line ${\mathbf c} = f(\chi)$
  with a slope $> 8.76$.  For example,
  $N$ is equal to $267145kx^2 + 70$  for some large
  numbers $k$ and $x$ determined by the condition
  $\frac{{\mathbf c}(X)}{\chi(X)} > 8.76$.\\

  Next we prove the following claim  by computing SW-basic classes of spin
  symplectic $4$-manifolds shown in Theorem~\ref{main-1}. \\

\hspace*{-1.8em} {\bf Claim 2}. {\em  None of the simply connected
 spin symplectic $4$-manifold in Theorem~\ref{main-1} admits
 a complex structure.} \\

  In order to compute the SW-basic classes of spin symplectic $4$-manifolds
 in Theorem~\ref{main-1}, we use the following two product formulas
 for the Seiberg-Witten invariant of a fiber sum $4$-manifold.

\begin{theorem}[\cite{mms}]
\label{mms}
  Suppose $X$ and $Y$ are smooth $4$-manifolds with \mbox{$b_{2}^+ >1$}
 which contain an embedded torus $f$ in a cusp neighborhood.
 Then the Seiberg-Witten invariant of a fiber sum $4$-manifold
 $X\sharp_{f}Y$ is given by
 \[ SW_{X\sharp_{f}Y} = SW_{X} \cdot SW_{Y} \cdot (e^{f} - e^{-f})^{2} . \]
\end{theorem}

\begin{theorem}[\cite{mst}]
\label{mst}
 Let $X$ and $Y$ be closed oriented smooth $4$-manifolds with
 $b_{2}^{+}(X) \geq 1$ and $b_{2}^{+}(Y) \geq 1$ which contain
 a smoothly embedded surface $\Sigma$ with genus $g >1$ representing
 a homology class of infinite order and of square zero.
 If there are characteristic classes  $l_{1}\in H^{2}(X;{\bf Z})$ and
 $l_{2}\in H^{2}(Y;{\bf Z})$ with $<l_{1},[\Sigma] > =<l_{2},[\Sigma] >=2g-2$
 and with $SW_{X}(l_{1}) \neq 0$ and $SW_{Y}(l_{2}) \neq 0$,
 then there exists a characteristic class
 $k\in H^{2}(X\sharp_{\Sigma} Y;{\bf Z})$ with
 $k|_{\partial D^{2} \times \Sigma}=\mathrm{proj}^{\ast}(k_{0})$ for
 $k_{0} \in  H^{2}(\Sigma ;{\bf Z})$ satisfying $<k_{0},[\Sigma ]> =2g-2$
 for which $SW_{X\sharp_{\Sigma} Y}(k) \neq 0$.
\end{theorem}

 \hspace*{-1.8em} {\em Remark.} \ Theorem~\ref{mms} tells us the complete
 Seiberg-Witten invariant of $X\sharp_{f}Y$ obtained by a taking a fiber
 sum along a torus $f$ in a cusp neighborhood. On the other hand,
 Theorem~\ref{mst} implies that there exists at least one SW-basic
 class of $X\sharp_{\Sigma}Y$ obtained by  taking a fiber sum along
 an embedded surface $\Sigma$ of genus $>1$ under an appropriate
 condition.\\

\hspace*{-1.8em} {\em Proof of Claim  2.} \  Note that spin
 symplectic $4$-manifolds $X'$, $Z$, $X$ and $W$ in the proof
 of Theorem~\ref{main-1} are of the forms
\begin{eqnarray*}
  Z &=&  X'\sharp_{\Sigma} X'  \, \, \, \mathrm{with}\, \,
    X' := X(2, 2g+1, 4g+1)\sharp \overline{\bf  CP}^{2}   \\
  X &= & \overbrace{Y\sharp_{\Sigma_{g}} \cdots \sharp_{\Sigma_{g}} Y}^{k}
         \sharp_{\Sigma_{g}}  Z     \\
  W &= & \overbrace{X\sharp_{f} X\sharp_{f} \cdots \sharp_{f} X}^{m} \sharp_{f}
       H(8k' -1) \sharp_{f} E(2n) \, \, \, \, \, \, \, \mathrm{or} \\
    &=& \overbrace{X\sharp_{f} X\sharp_{f} \cdots
      \sharp_{f} X}^{m} \sharp_{f} H(7)\sharp_{f} H(8k' -1) \sharp_{f} E(2n)
\end{eqnarray*}

\hspace*{-1.75em} for some positive integers $g, m, n, k'$ and $k \in {\mathbf Z}$.

 Now let us try to find SW-basic classes of $X'$, $Z$, $X$, and $W$
 consecutively: First, since the manifold $X'$ is a blow up of an
 algebraic surface $X(2, 2g+1, 4g+1)$ of general type and
 $X(2, 2g+1, 4g+1)$ has only one up to sign SW-basic class $K$,
 the canonical class of $X(2, 2g+1, 4g+1)$, the SW-basic classes
 of $X'$ are up to sign of the forms $\{K + {\mathbf e},
 K -{\mathbf e}\}$, where ${\bf e}$ is a class represented by
 an exceptional curve in $\overline{{\bf CP}}^{2}$.
 Secondly,  since  $Z$ is the $4$-manifold obtained by taking a fiber sum
 along $\Sigma = T-{\mathbf e}$ and $\Sigma$ is a symplectically embedded
 surface in $X'$, it satisfies the hypotheses of Theorem~\ref{mst} above.
 Hence, by applying Theorem~\ref{mst} or V. Mu$\mathrm{\tilde{n}}$oz' result
 in~\cite{mu}, it is easy to check that only a pair $( K + {\bf e}, K + {\bf e})$
 of SW-basic classes of $(X',X')$ induces the SW-basic class of $Z$,
 denoted by  $K_{Z}$.
  Next, since $Y$ is also a minimal complex surface of general type,
 it has only one SW-basic class up to sign, say $K_{Y}$. Then,
 by applying Theorem~\ref{mst} on a fiber sum manifold $Y\sharp_{\Sigma_{g}} Z$
 again, we obtain one SW-basic class of $Y \sharp_{\Sigma_{g}} Z$ induced from
 $(K_{Y}, K_{Z})$, which is denoted by $K_{Y \sharp_{\Sigma_{g}} Z}$.
 By repeating the same procedures as above,
 we get one SW-basic class up to sign, say $K_{X}$, of a fiber sum manifold $X$.
  Finally, let us compute SW-basic class of our manifold $W$. Since Horikawa
 surfaces $H(7)$ and $H(8k'-1)$ have only one SW-basic class, say $K_{H_{7}}$ and
 $K_{H}$, respectively, and since an elliptic surface $E(2n)$ has the following
 SW-basic classes $\{ 2kT \, | \,  k= 0, 1, \, \ldots, \, (n-1) \}$
 where $T$ is a regular torus fiber in the nucleus $C(2n) \subset E(2n)$,
 by applying Theorem~\ref{mms} repeatedly on a Lagrangian torus $f$
 embedded in Brieskorn submanifolds of $X$, $H(7)$, $H(8k'-1)$ and $Z$,
 we obtain SW-basic classes up to sign of $W$:
 \[ K_{H} \pm({K_{H_{7}}}) \pm{2k}T \pm
    \overbrace{K_{X} \pm \cdots \pm{K_{X}}}^{m}
    + \overbrace{\{0,\pm{2}f\}+\cdots +\{0,\pm{2}f\}}^{m+1 \, \,
     (\mathrm{or,} \, \, m+2)}  . \]
 Note that  all these manifolds, $W$, have obviously  more than one SW-basic class.
 This means that any such a manifold $W$ with $\mathbf{c}(W) > 0$ cannot be
 diffeomorphic to a complex surface of general type. In the case when
 ${\mathbf c}= 0$, we already constructed a family of simply connected spin
 non-complex symplectic $4$-manifolds $X(2n)$ covering the elliptic line in
 the proof of Proposition~\ref{main-pro}. Hence Claim 2 is proved. $\ \ \ \Box$
\\

\hspace*{-1.8em} {\em Acknowledgement}.
 The author would like to thank Ronald Fintushel and Andr\'{a}s Stipsicz
 for helpful conversations. He is indebted to Ronald Fintushel for a critical
 reading of the first draft. He would also like to thank Chris Peters for
 explanation on spin complex surfaces they constructed in~\cite{ppx}.
 Finally, he also wishes to thank the referee for pointing out some errors
 in the earlier version of this paper.\\
 \\

\end{document}